# Artificial Neural Networks in Fluid Dynamics: A Novel Approach to the Navier-Stokes Equations

Extended Abstract

Megan F. McCracken
Austin Peay State University
Department of Physics, Engineering and Astronomy
601 College St., Clarksville, TN 37044
mmccracken1@my.apsu.edu

## ABSTRACT

Neural networks have been used to solve different types of large data related problems in many different fields. This project takes a novel approach to solving the Navier-Stokes Equations for turbulence by training a neural network using Bayesian Cluster and SOM neighbor weighting to map ionospheric velocity fields based on 3-dimensional inputs. Parameters used in this problem included the velocity, Reynold's number, Prandtl number, and temperature. In this project data was obtained from Johns-Hopkins University to train the neural network using MATLAB. The neural network was able to map the velocity fields within a 67% accuracy of the validation data used. Further studies will focus on higher accuracy and solving further non-linear differential equations using convolutional neural networks.

## CCS CONCEPTS

• **Artificial Neural Networks → SOM neighbor-weighting model**; *Fluid Dynamics*; Velocity Fields

## KEYWORDS

ACM proceedings, Fluid Dynamics, Velocity Fields, Neural Networks, Navier-Stokes

## 1 INTRODUCTION

Modeling the dynamics of non-linear fluid flow has been difficult to achieve with current technologies. A small portion of research in fluid dynamics has been devoted to modeling different types of fluid flows through computer simulations. This is known as CFD, or Computational Fluid Dynamics. Utilizing the ability of computers to achieve more efficient computation, simulations of annealing, civil engineering, and weather predictions have been created using Computational Fluid Dynamics [1]. CFD has grown in its ability to accurately model fluid flow and predict paths of the fluid being modeled, however, the complexity of the systems being observed bely a certain amount of inaccuracy within the models. Recently researchers in the field of fluid dynamics have been looking at neural networks and their ability to solve complex problems quickly. The most commonly used model to date has been the Feed Forward Neural Network [2]. This model takes several inputs then feeds the information forward through hidden layers and produces an output. In order to optimize the output, the network uses a technique known as backpropagation. This technique changes the weights of certain "neural pathways" to make them impactful to the final output. Neural Networks have been used sparingly in CFD, because of the non-linear nature of the Navier-Stokes equations. This complexity obviates the use of basic feed-forward neural networks within CFD. Therefore, there have been strides taken to increase the complexity of space in which neural networks can work [3], including Bayesian Cluster Neural Networks, and Self-Organized Mapping (SOM) Neural Networks, which were applied in this research.

## 2 EXPERIMENTAL AND COMPUTATIONAL DETAILS

### 2.1 Bayesian Cluster Neural Networks

Feed Forward Neural Networks work by using training data sets fed through hidden layers, then backpropagated to reweight. After the network has been trained, a validation set will be used as proof of the accuracy and test data sets can be run through the model. This process allows the Neural Network to not overfit the data it was initially given. Bayesian Cluster Neural Networks work in a similar fashion, but instead of having one connection between nodes, the weights can be changed in between the hidden layers. This allows the Bayesian Cluster Neural Networks to change in a non-linear fashion, according to variance in the system. This model takes a longer time for the computer to process but produces a more accurate result. Both the Bayesian Cluster method and the SOM neighbor-weighting model utilize this cluster mapping method to format a result.

## 2.2 SOM Neighbor-Weighting

The final product produced for the research took a two-tiered approach. First, the Bayesian Cluster model was implemented then to solidify the accuracy of the output SOM feature mapping was applied.

Self-Organized Map Neighbour-weighting, or Feature Mapping, namely the Kohonen method [4] was used to analyse the final product produced by the data sets, after using the Bayesian Cluster method. This method weights nodes based on the likelihood that they are close to the predicted value, and produces features based on a large amount of 'hits' for that specific feature.

This is a common technique used in facial recognition, because it produces peak nodes on prominent features such as a nose, or high cheekbones. For this research the neural network was optimized to recognize the highest Reynold's number, or highest amount of 'turbulence' within the system. The results of which can be seen in the results section of this paper.

## 2.3 Data Set Acquisition

Ionospheric data was taken from the Johns-Hopkins Atmospheric Data Center [5] and processed via MATLAB [6]. The data was extracted at fixed points in time, which allowed for the data to be analyzed without the time component. The Navier-Stokes Equations become significantly more complex when solved over time and were not within the scope of the research.

CFD modelling focuses on the ability to predict a fluid flow through time and utilizes the Navier-Stokes Equations to better understand turbulence within the system [7]. The data acquired from Johns-Hopkins University was a snapshot of the fluid flow within one specific section of time and did not consider the flow between the times the data was accessioned. This is not to say that the model produced in this research would not be able to eventually model the fluid flow through time, but that the predictions produced are not time dependent, and are to be considered instantaneous.

## 2.4 The Navier-Stokes Equations

*2.4.1 Format.* The Navier-Stokes Equations, shown in Fig. 1, are incredibly complex, and not well understood. The basic format of the Navier-Stokes equation utilizes the several different variables to produce a model of fluid flow with the dependence on turbulence which was not included in the previous Euler models [1].

$$\begin{aligned}
\frac{\partial u}{\partial t} + u\frac{\partial u}{\partial x} + v\frac{\partial u}{\partial y} + w\frac{\partial u}{\partial z} &= -\frac{\partial P}{\partial x} + 1/Re \left(\frac{\partial^2 u}{\partial x^2} + \frac{\partial^2 u}{\partial y^2} + \frac{\partial^2 u}{\partial z^2}\right), \\
\frac{\partial v}{\partial t} + u\frac{\partial v}{\partial x} + v\frac{\partial v}{\partial y} + w\frac{\partial v}{\partial z} &= -\frac{\partial P}{\partial y} + 1/Re \left(\frac{\partial^2 v}{\partial x^2} + \frac{\partial^2 v}{\partial y^2} + \frac{\partial^2 v}{\partial z^2}\right), \\
\frac{\partial w}{\partial t} + u\frac{\partial w}{\partial x} + v\frac{\partial w}{\partial y} + w\frac{\partial w}{\partial z} &= -\frac{\partial P}{\partial z} + 1/Re \left(\frac{\partial^2 w}{\partial x^2} + \frac{\partial^2 w}{\partial y^2} + \frac{\partial^2 w}{\partial z^2}\right), \\
\frac{\partial u}{\partial x} + \frac{\partial v}{\partial y} + \frac{\partial w}{\partial z} &= 0.
\end{aligned}$$

**Figure 1: Navier-Stokes Equations in the cartesian coordinate system, as well as the continuity equation**

The variables used in the Navier-Stokes equations include Reynold's number, Prandtl number, velocity in each direction, pressure, temperature, and the dependence of each of these through time [8]. When applying the Neural network, each of these variables was an input and was weighted across 10 hidden layers.

*2.4.2 Managing the Input Values.* The number of values taken was large enough to optimize the data within the means of the Navier-Stokes equations. Because there is a high dependence of each parameter in the Navier-Stokes Equation to all other parameters, the optimization of the output was best achieved when culled to a smaller data set manageable for a single processor. Further research will be focused on larger data sets and parameterization. When utilizing the SOM model, the data was specifically parameterized to focus on high-Reynold's number areas allowing the computer to be able to process the large number of data points taken from each direction within the velocity fields (Fig. 2).

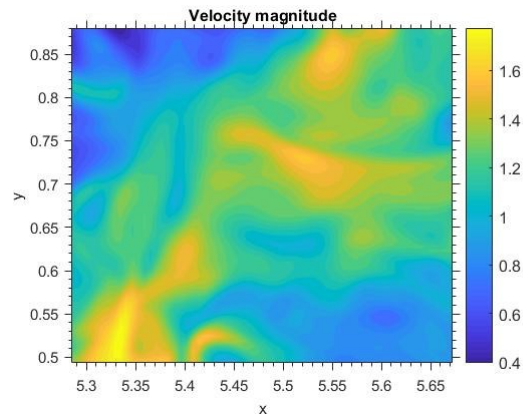

**Figure 2: A velocity field of one of the test data sets used to train the neural networks.**

## 3 RESULTS AND DISCUSSION

### 3.1 Bayesian Cluster Output

The results for the Bayesian Cluster Model were inconclusive based on the 496 epochs out of 500 that were run. Both models utilized the method of 70% training, 15% validation, and 15% test. This means that out of the data used, which included 4096 inputs, 2867 were used to train the model and 614 were used in validation and testing. The model produced a steady state until the final epochs where the mu values increased severely and could no longer continue running through the model. This result can be seen in Fig. 3 and 4 below, which documents the output of each of the epochs according to the iterations run.

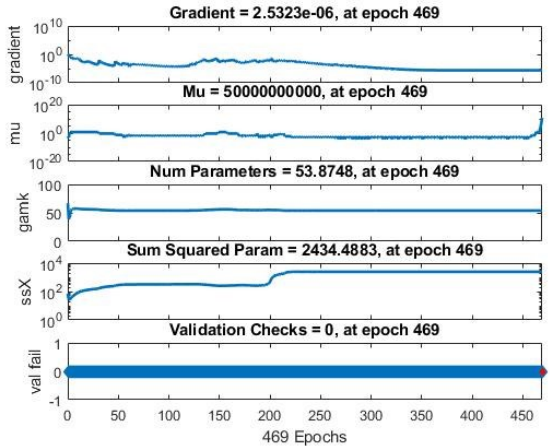

**Figure 3: The steady state of the Bayesian Cluster model. Note, at the tail end of the mu plot the value increases drastically. This caused the neural network to not be able to produce a viable output.**

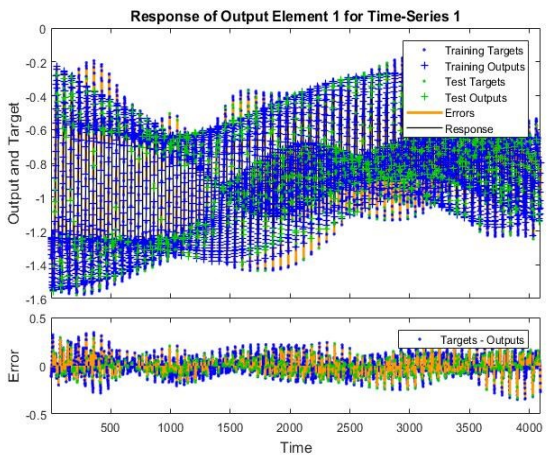

**Figure 4: Output elements based on the training and test data sets.**

Fig. 4 shows the output for the training and test sets within the Bayesian Cluster model. Given the large number of data points and the high amount of parameterization that the neural network was attempting to achieve, this model was not able to validate the output, and/or over-parameterized the model and could not find an optimal state. The network was able to cluster the outputs initially, allowing them to be more easily used in the second tier of the neural network model.

### 3.2 SOM Output

The data taken from the Bayesian Cluster model was then used as the initial input of the SOM portion of the neural network. This was an initial novel approach to solve the over-parameterization of the Bayesian Cluster output set and recreate velocity fields utilizing the Reynold's number as the strongest feature.

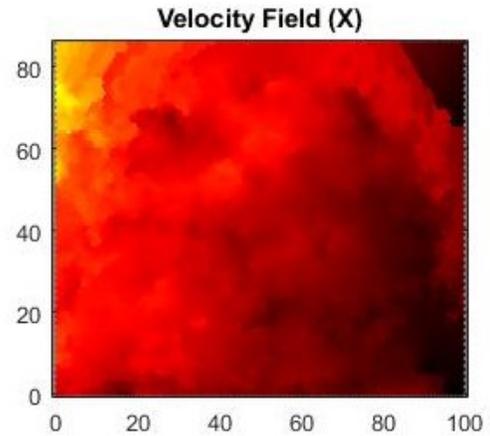

**Figure 5: SOM Feature map produced for the X direction.**

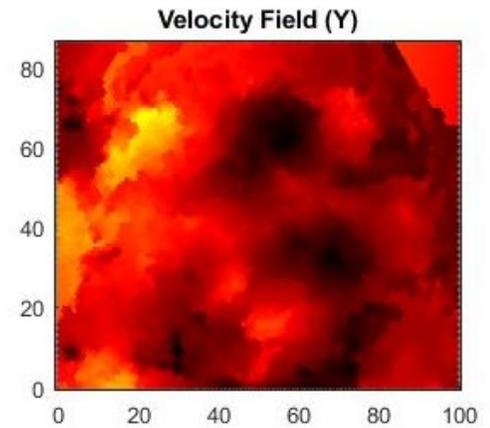

**Figure 6: SOM Feature map produced for the Y direction.**

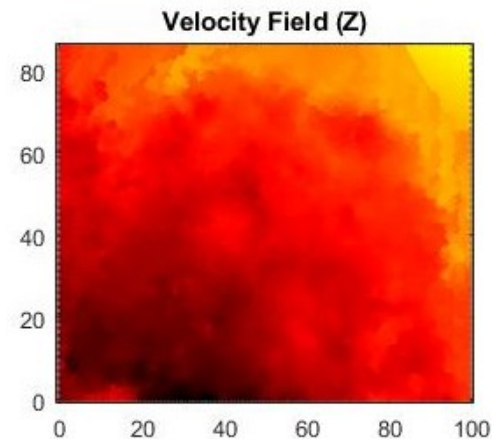

**Figure 7: SOM Feature map produced for the Z direction.**

The velocity Field outputs for the SOM neighbor-weighting model were not easily concatenable and produced in three different figures to show the correlation between the produced velocity fields and the test data set (Fig. 8). The size of the data sets hindered the ability for the model to be tested multiple times, which is further discussed in the conclusions section, however the two-tier model managed to achieve an accuracy of .67, to that of the validation sets.

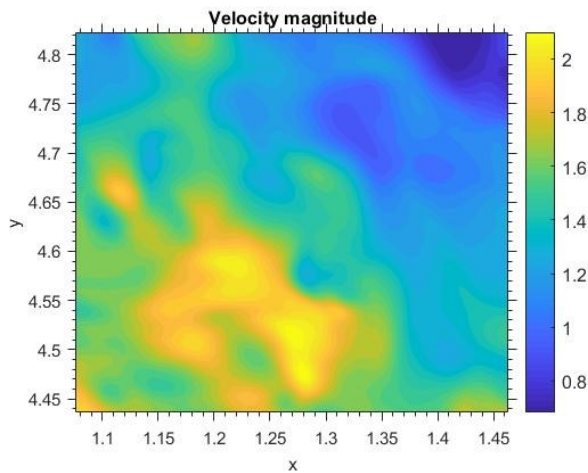

**Figure 8: Visual representation of the data used to validate the SOM neighbor-weighted model.**

## 3.3   Analysis of Model Accuracy

As discussed previously, this research took a unique approach to modeling the flow of a fluid using the Navier-Stokes equations. The application of the two different models of neural network allowed for a higher accuracy rate to be achieved using the Reynold's number as the most important feature. Other parameters could have been used as the main feature used in the SOM model and would achieve different results depending on the parameter used. The two-tiered approach used in this research achieved an accuracy rate of .67, using the six training sets fed through the Neural Network. Improvements that could be made to this research include a larger amount of validation sets, parallel processing and validation through time. However, due to a time and resource limitation this research was only validated using one set of data. To achieve a higher accuracy and trust in the model output the two-tier set up should be run multiple times, while achieving similar results.

## 4   CONCLUSIONS

In summary, the ability for neural networks to solve the reproduce an accurate model of fluid flow within a finite space was reproduced using a two-tier approach to Neural Networks that has not previously been attempted. Computational Fluid Dynamics researchers have tried to accurately recreate fluid flow over time using different models and have not been able to achieve a reproducible turbulence model. The Bayesian Cluster and the SOM feature neural network models used in tandem were able to reproduce velocity fields within a 67% accuracy of the actual output. This model of the Navier-Stokes equations was not the most accurate model produced it can be expanded upon to produce a higher optimization, even given the complexity of the problem. This model can be expanded to other non-linear differential equation models or can be organized to attain a higher accuracy within this model itself. The research presented here took a novel approach to solving a complex non-linear problem and was able to produce accurate feature maps. The model will be expounded upon to produce a higher accuracy, with the utilization of more processors and higher parameterization.

## ACKNOWLEDGMENTS

This research was made possible by Johns-Hopkins University Computational Fluid Dynamics data sets, without which the models could not have been created. Also, a special thanks to Dr. Daniel Mayo and Dr. Justin Oelgoetz who helped in the explanation and understanding of neural networks.